\documentclass[12pt]{amsart}

\usepackage[textheight=23.5cm, left=1.2cm, right=1.2cm]{geometry}

\usepackage[english]{babel}
\usepackage[T1]{fontenc}
\usepackage[utf8]{inputenc}

\theoremstyle{plain}
    \newtheorem{theorem}{Theorem}[section]
    \newtheorem{lemma}[theorem]{Lemma}

    \newtheorem{conjecture}[theorem]{Conjecture}

\theoremstyle{definition}

\usepackage[
    draft = false,
    unicode = true,
    colorlinks = true,
    allcolors = blue,
    hyperfootnotes = true,
    citecolor = red
]{hyperref}
\usepackage{amsfonts,amsmath,amssymb,amscd,amsthm,url}
\usepackage[inline]{enumitem}
\usepackage{graphicx,epsf,afterpage}
\usepackage{soul}
\usepackage{multicol}
\usepackage{array}
\usepackage{braket}
\usepackage{epigraph}
\usepackage{rotating}

\usepackage[
backend=biber,
style=numeric-comp,
sorting=none,
giveninits=true
]{biblatex}
\usepackage{csquotes}

\usepackage{xcolor}

\usepackage{ulem}

\usepackage{tikz}
\usetikzlibrary{positioning, calc, intersections, through, arrows, matrix, chains, math}
\usetikzlibrary{shadows}
\usetikzlibrary{graphs}
\usetikzlibrary{graphs.standard}
\usetikzlibrary{patterns}
\usetikzlibrary{decorations.pathreplacing,calligraphy,backgrounds}

\newcommand\norm[1]{\ensuremath{\left\lVert#1\right\rVert}}

\renewcommand{\Pr}{\mathrm{P}}

\DeclareMathOperator{\Expect}{\mathbb{E}}

\DeclareMathOperator{\var}{var}

\newcommand{\Lcal}{\mathcal{L}}
\newcommand{\Hcal}{\mathcal{H}}
\newcommand{\Tcal}{\mathcal{T}}

\newcommand{\R}{\ensuremath{\mathbb{R}}}

\newcommand{\Cx}{\ensuremath{\mathbb{C}}}

\newcommand{\N}{\ensuremath{\mathbb{N}}}

\renewcommand{\geq}{\geqslant}

\renewcommand{\leq}{\leqslant}

\newcounter{mcnt}

\newcounter{wordcnt}

\begin{document}

\title[The law of large numbers for discrete generalized quantum channels]{The law of large numbers \\ for discrete generalized quantum channels}

\author{S.\,V. Dzhenzher and V.\,Zh. Sakbaev}

\begin{abstract}
     We consider random linear operators $\Omega \to \mathcal{L}(\Tcal_p, \Tcal_p)$ acting in a $p$-th Schatten class $\Tcal_p$ in a separable Hilbert space $\Hcal$ 
     for some $1 \leqslant p < \infty$. Such a superoperator is called a pre-channel since it is an extension of a quantum channel to a wider class of operators without requirements of trace-preserving and positivity.
     Instead of the sum of i.i.d. variables there may be considered the
     composition of random semigroups $e^{A_i t/n}$ in the Banach space $\Tcal_p$.
     The law of large numbers is known in the case $p=2$ in the form of the usual law of large numbers for random operators in a Hilbert space.
     We obtain the law of large numbers for the case $1\leq p \leqslant 2$.
\end{abstract}

\maketitle
\thispagestyle{empty}

\tableofcontents

\section{Introduction}

The theory of products of random linear operators is constructed and developed in a lot of publications \cite{Meta, Furstenberg, Tutubalin}. 
The theory of random  variables with values in a space of linear
operators in Banach spaces, and random processes with values in the
Banach space of bounded linear operators has applications in different problems of mathematical physics.
A random quantum dynamics can be considered either as a
composition of random linear operators in a Hilbert space of a quantum
system or as a composition of random quantum channels in the set of quantum states \cite{Aaron, Hol}.
A random operator valued processes arise in the study of the ambiguity of the quantization procedure
for Hamiltonian systems \cite{Berezin, OSS-2016, Orlov2}, of quantum random walks \cite{Aaron, Hol, Pechen-Volovich}, and
of dynamics of quantum systems interacting with a quantum statistical
ensemble \cite{GOSS21, GOSS-2022, Joye, Pechen}.

Distributions of
products of random linear operators had been studied in
\cite{Oseledets, Tutubalin, Skorokhod, New-23}. The law of large
numbers (LLN) for a composition of independent identically distributed (i.i.d.) random linear operators in
a Banach space had been formulated in terms of distributions of spectral
characteristics of compositions.
The analogue of the LLN which is formulated for compositions of
independent random linear operators is studied as the probability estimation for the deviation of a
random composition from its mean value in different operator topologies \cite{S16, S18, OSS19}. The LLN for operators composition is the non-commutative generalization of the LLN for the sum of independent random vector valued variables. Different versions of limit theorems of LLN and CLT type for compositions of i.i.d. random operators are obtained in \cite{Berger, GOSS-2022, SSSh23}.

A lot of results on the LLN \cite{S16, S18, OSS19} are obtained
for random linear operators acting in a Hilbert space. But a
consideration of quantum dynamics of an open quantum system needs to study the random mappings in the Banach space of quantum states \cite{Hol, GOSS-2022, KOS, Pechen, Sahu}. However, the case of random operators
acting in a Banach space is studied significantly less. Some results
were obtained on the LLN for semigroups of bounded linear operators acting in a finite dimensional space \cite{GOSS-2022}.
Random strongly continuous semigroups of bounded linear operators acting in Banach spaces of sequences $\ell_p,\ p\geqslant 1$, are considered in \cite{DzhenzherSakbaev24}. The LLN had been obtained for  compositions of random semigroups of linear operators acting in the spaces \(\ell_p, \ p\in [\,1,2\,]\). In this paper we study a non-commutative extension of paper \cite{DzhenzherSakbaev24} results for  compositions of random semigroups of linear operators acting in the Schatten classes \(\Tcal_p, \ p\in [\,1,2\,]\).


The studying of random transforms of a Hilbert space of a quantum
system leads to the analysis of compositions of random unitary
operators acting in the Hilbert space. The theory of random
transformation of a Hilbert space is developed in \cite{Kempe,
OSS19}.
The interest for composition of random linear operators in a Banach
space arises in studying of random quantum channels \cite{GOSS-2022,
New-23, N-23}. To this aim we should consider random variables in the set of quantum channels. 
A quantum channel is a completely positive trace-preserving linear mapping acting in a Banach space of nuclear operators. 

In this paper we study random generalized quantum channels on discrete probability space; we call them random pre-channels for short.
So, a pre-channel is a bounded linear operator $A\in \Lcal(\Tcal_p)$ acting in a Schatten class $\Tcal_p$ with some $p\in [\,1,+\infty )$. 
We obtain the LLN for compositions of i.i.d. random pre-channels in the Banach space $\Tcal_p$ for $p\in [\,1,2\,]$.

The crucial role in the obtaining of LLN is the proof of Chebyshev's estimate for the probability of the deviation of products of random pre-channels from their expected values in a suitable topology. 
To solve the problem of obtaining Chebyshev's estimates, we should prove the measurability property for the product of two random pre-channels and for the conjugate channel acting in the space $\Tcal_q$ where $q=\frac{p}{p-1}$. To avoid these problems of measurability, we consider pre-channel valued random variables with a discrete probability space. The topology of the LLN validity for a sequence of products  of i.i.d. random pre-channels in the Banach space $\Tcal_p$ with $p\in [\,1,2\,]$ is chosen as the topology on the space $\Tcal_p$ which is generated by the collection of functionals that generate the strong operator topology (SOT) of the space $\Lcal(\Tcal_2)$. Thereby, the inclusion $\Tcal_p\subset \Tcal_2,\ p\in [\,1,2\,]$ is a reason to consider the LLN for compositions of random pre-channels in the Banach spaces $\Tcal_p$ with $p\in [\,1,2\,]$ only.

The structure of the article is the following.
In~\S\ref{s:main-res}, the
main result (Theorem~\ref{t:wlln}) is given.
In~\S\ref{s:lemmas}, some auxiliary lemmas are given.
In~\S\ref{s:proof}, the proof of Theorem~\ref{t:wlln} is given.

\section{Background and statements}\label{s:main-res}

\subsection*{Topology}

In this text we consider Banach spaces over some field $\mathbb{K}$ ($\R$ or $\Cx$).
Recall that for Banach spaces $X,Y$ the space $\Lcal(X,Y)$ is the Banach space of linear bounded operators $X \to Y$.
We shortly denote \(\Lcal(X) := \Lcal(X,X)\).

Let \(\Hcal\) be a separable Hilbert space.
For \(1 \leq p < \infty\) let \(\Tcal_p := \Tcal_p(\Hcal) \subset\Lcal(\Hcal)\) be the Banach space of $p$-th Schatten-class operators.
For convenience, we identify \(\Tcal_\infty\) with the space of bounded operators.
We denote norms in these spaces simply as \(\norm{\cdot}_p\).

For avoiding collisions, we will call the operators from \(\Lcal(\Tcal_p, \Tcal_q)\) as pre-channels.
Note that in quantum mechanics quantum channels are operators from \(\Lcal(\Tcal_1)\) with some additional requirements, which we do not need.
Talking about pre-channels, by default we mean that \(1\leq p,q\leq \infty\).
We denote norms in this space simply as \(\norm{\cdot}\).

\subsection*{Probability}

Let \((\Omega, \mathcal{F}, \Pr)\) be a discrete probability space.
Note that discreteness guarantees that all random objects we work with are measurable.

Let $X$ be a Banach space. Let $\xi\colon\Omega\to X$ be a random element (a pre-channel or a $p$-th Schatten-class operator in our cases).
Its \textbf{\emph{expected value}}
\[
    \Expect\xi := \sum_{x \in X} x \Pr\{\xi=x\}
\]
is the Bochner integral, where the sum converges almost surely in the norm topology of $X$. It is well-known \cite[Theorem~2]{DiestelUhl} that $\xi$ is Bochner-integrable if and only if \(\norm{\xi}\) is integrable.

\begin{lemma}[Integration]
    Let \(A \colon \Omega \to \Lcal(\Tcal_p,\Tcal_q)\) be a random pre-channel with integrable norm.
    Then for any \(x\in\Tcal_p\)
    \[
        (\Expect A)x = \Expect (Ax).
    \]
\end{lemma}
\begin{proof}
    We have
    \(\displaystyle
        (\Expect A)x = \sum_{A_n} A_nx\Pr\{A=A_n\} = \sum_y y\Pr\{Ax=y\} = \Expect (Ax).
    \)
\end{proof}

We say that random channels \(\{A_n\}\) are \textbf{\emph{independent}} if for any corresponding
non-random pre-channels \(\{B_n\}\)
and any chosen \(i_1,\ldots,i_k\) we have
\[
    \Pr\{A_{i_1} = B_{i_1},\,\ldots,\,A_{i_k} = B_{i_k}\} = \Pr\{A_{i_1} = B_{i_1}\}\ldots\Pr\{A_{i_k} = B_{i_k}\}.
\]

\subsection*{Semigroups}

For a bounded pre-channel \(A\colon\Omega\to\Lcal(\Tcal_p)\) denote by \(e^{At}\) the random (uniformly continuous) semigroup with $A$ as the infinitesimal generator.

A sequence of compositions $W_n(t)$ of random semigroups generated by a random pre-channels $\Omega\to\Lcal(\Tcal_p)$
\emph{converges almost surely in SOT of $\Tcal_q$ to a random semigroup $W(t)$ uniformly for $t > 0$ in any segment} if for any $x \in \Tcal_p$ and $T > 0$
\[
    \lim\limits_{n\to\infty} \sup\limits_{t \in [\,0,T\,]} \norm{({W_n(t) - W(t)})x}_q = 0
    \quad\text{a.s.}
\]

The following theorem is known in the case of unitary semigroups generated by self-adjoint operators. We consider the case of bounded generators.

\begin{theorem}[Main result; proved in \S\ref{s:proof}]\label{t:wlln}
    Let $1 \leq p \leq 2$.
    Let $A\colon\Omega\to\Lcal(\Tcal_p)$ be a generator of a semigroup $e^{At}$ such that $\norm{A}$ is bounded.
    Let $A,A_1, A_2, \ldots$ be a sequence of random i.i.d. generators.
    Then $e^{A_1t/n} \ldots e^{A_nt/n}$ converges almost surely in SOT of $\Tcal_2$ to $e^{\Expect At}$ uniformly for $t$ in any segment.
\end{theorem}

It is interesting to obtain the LLN in SOT of $\Tcal_p$ not $\Tcal_2$.
This leads us to the following conjecture.

\begin{conjecture}
    Under the assumptions of Theorem~\ref{t:wlln}, $e^{A_1t/n} \ldots e^{A_nt/n}$ converges almost surely in SOT of $\Tcal_p$ to $e^{\Expect At}$ uniformly for $t$ in any segment.
\end{conjecture}

\section{Auxiliary lemmas}\label{s:lemmas}

For a Banach space $X$ we denote by $X^*$ the \emph{dual space} of linear functionals \(X\to\mathbb{K}\).
For \(x \in X\) and \(f \in X^*\) we identify \(\braket{f, x} := f(x)\).

For any \(1 \leq p \leq \infty\) we define \(p^*\) by the identity \(\frac{1}{p} + \frac{1}{p^*} = 1\).
So, for \(p < \infty\) we have \(\Tcal_p^* = \Tcal_{p^*}\).

For \(U \in \Lcal(\Tcal_p, \Tcal_q)\) define the \emph{adjoint pre-channel} \(U^* \in \Lcal(\Tcal_q^*, \Tcal_p^*)\) by
\[
    \braket{U^*y, x} = \braket{y, Ux}
\]
for any \(x\in \Tcal_p\), \(y\in\Tcal_q^*\).

\begin{lemma}[Adjoint]\label{l:adjoint}
    Let \(1\leq p < \infty\).
    Let \(A\colon\Omega\to\Lcal(\Tcal_p)\) be a random pre-channel with integrable norm.
    Then its adjoint $A^*$ is integrable, and
    \[
        \Expect A^* = (\Expect A)^*.
    \]
\end{lemma}
\begin{proof}
    The adjoint $A^*$ is integrable since \(\norm{A^*}=\norm{A}\) is integrable.
    Take any \(x\in\Tcal_p\) and \(y\in\Tcal_p^*\).
    Consider the linear functional \(T_x\in\Tcal_{p^*}^*\) defined by \(T_xf := \braket{f,x}\).
    Since \(\Expect(A^*y)\) is the Bochner integral, we have that
    \(T_x (A^*y)\) is also Bochner integrable, and
    \[
        \Expect T_x(A^*y) = T_x \Expect (A^*y).
    \]
    Then
    \[
        \braket{(\Expect A^*)y, x} = \Expect\braket{A^*y,x} = \Expect\braket{y, Ax} = \braket{y, (\Expect A)x} = \braket{(\Expect A)^*y,x},
    \]
    where the equality before the last follows analogously (or since the expected value is the Pettis integral).
\end{proof}

\begin{lemma}[Super-operator]\label{l:super-op}
    Let \(1\leq p\leq 2\).
    Let \(A,B\colon\Omega\to\Lcal(\Tcal_p)\) be independent random pre-channels with bounded norms.
    Let \(\Expect A = 0\).
    Then
    \[
        \Expect B^*AB = 0.
    \]
\end{lemma}
\begin{proof}
    For a pre-channel \(K\in\Lcal(\Tcal_p)\) consider a super-operator \(U_K \colon\Lcal(\Tcal_p)\to\Lcal(\Tcal_p,\Tcal_p^*)\) defined by
    \[
        U_K T := K^*TK.
    \]
    The random pre-channel \(U_BA\) has the bounded norm and hence integrable.
    Then 
    \[
        \Expect (U_BA) = \sum_{K,T} U_KT \Pr\{B=K\}\Pr\{A=T\} =
        \sum_K \Pr\{B=K\} U_K \sum_T T\Pr\{A=T\} = (\Expect U_B) (\Expect A) = 0.
    \]
\end{proof}

\begin{lemma}[Independence]\label{l:indep}
    Let \(A_1, \ldots, A_n\) be independent random pre-channels with bounded norms such that their composition is well-defined (so that their images and domains are correspondent).
    Then
    \[
        \Expect (A_1\ldots A_n) = (\Expect A_1)\ldots (\Expect A_n).
    \]
\end{lemma}
\begin{proof}
    Since for independent random pre-channels compositions of different operators are also independent, it is sufficient to prove the lemma for $n=2$.
    For two random pre-channels $A,B$
    \begin{multline*}
        \Expect AB = \sum_T T \Pr\{AB = T\} =
        \sum_{T,U} T \Pr\{B = U\}\Pr\{AU=T\} =
        \sum_U \Pr\{B = U\} \Expect (AU) = \\ =
        (\Expect A) \sum_U U \Pr\{B = U\} = (\Expect A)(\Expect B).
    \end{multline*}
\end{proof}

Note that the proofs of the last two lemmas are similar but not the same, since in Lemma~\ref{l:super-op} we have the action of the super-operator $U_B$, while in Lemma~\ref{l:indep} we have the composition of two pre-channels.

\section{Proof of Theorem~\ref{t:wlln}}\label{s:proof}

Recall that $\Tcal_p \subset \Tcal_q$ for $p < q$.

Let $1 \leq p \leq 2$.
For a random pre-channel $A\colon\Omega\to \Lcal(\Tcal_p)$ with bounded norm denote
\[
    \var A := \Expect \bigl((A - \Expect A)^* (A - \Expect A)\bigr) \in \Lcal(\Tcal_p, \Tcal_p^*).
\]

\begin{lemma}[Chebyshev's inequality]
\label{l:cheb-l2}
    Let \(1 \leq p \leq 2\).
    Let \(A\colon\Omega\to \Lcal(\Tcal_p)\) be a random pre-channel with bounded norm.
    Then for any $\varepsilon>0$ and $x \in \Tcal_p$
    \[
        \Pr\bigl\{\norm{(A - \Expect A)x}_2 > \varepsilon\bigr\} \leq \frac{\norm{(\var A)x}_{p^*} \norm{x}_p}{\varepsilon^2}.
    \]
\end{lemma}
This lemma is known in the case $p=2$.
Its proof is analogous to the proof of \cite[Lemma~4 (Chebyshev's inequality)]{DzhenzherSakbaev24}, where it is proved for \(\ell_p\) instead of \(\Tcal_p\).

\begin{proof}[Proof of Theorem~\ref{t:wlln}]
    The idea of the proof is due to \cite[proof of Theorem~2]{GOSS-2022}.
    We will follow \cite[proof of Theorem~1]{DzhenzherSakbaev24}.
    
    Denote
    \[
        W_n(t) := e^{A_1t/n} \ldots e^{A_nt/n}.
    \]
    Fix $1 \leq p \leq 2$ and any $\varepsilon > 0$, $T > 0$ and $x \in \Tcal_p$.
    Since the probability space is discrete, it is sufficient to prove that
    \[
        \lim\limits_{n\to\infty} \Pr \left\{\sup\limits_{t \in [\,0,T\,]} \norm{\left(W_n(t) - e^{\Expect At}\right)x}_2 > \varepsilon\right\} = 0.
    \]
    The idea of the proof is to show the following two asymptotic equivalences in probability, which imply the theorem:
    \begin{align}
         \label{eq:lln}
         \lim\limits_{n\to\infty} \Pr \left\{\sup\limits_{t \in [\,0,T\,]} \norm{\left(W_n(t) - \Expect W_n(t)\right)x}_2 > \varepsilon\right\} &= 0, \\
         \label{eq:chernov}
         \lim\limits_{n\to\infty} \Pr \left\{\sup\limits_{t \in [\,0,T\,]} \norm{\left(e^{\Expect At} - \Expect W_n(t)\right)x}_2 > \varepsilon\right\} &= 0.
    \end{align}

    Equality~\eqref{eq:chernov} can be proved analogously to \cite[equality~(2)]{DzhenzherSakbaev24}.
    So we need to prove~\eqref{eq:lln}.
    For $n \in \N$ and $t \geq 0$ denote
    \[
        \Delta_n(t) := e^{A_nt} - \Expect e^{At}.
    \]
    For integer $0 \leq k \leq n$ and $1 \leq a_1 < \ldots < a_k \leq n$ denote
    \[
        F_{n;k;a_1, \ldots, a_k}(t) := (\Expect e^{At})^{a_1-1} \Delta_{a_1}(t) (\Expect e^{At})^{a_2-a_1-1} \ldots \Delta_{a_k}(t) (\Expect e^{At})^{n-a_k}.
    \]
    Then
    \[
        W_n(t) = \sum_{k=0}^n\sum_{1 \leq a_1 < \ldots < a_k \leq n} F_{n;k;a_1, \ldots, a_k}(t/n).
    \]
    Then
    \begin{multline*}
        \Expect W_n(t)^*W_n(t) =
        \sum_{k=0}^n\sum_{m = 0}^n\sum_{\substack{1 \leq a_1 < \ldots < a_k \leq n \\ 1 \leq b_1 < \ldots < b_m \leq n}} \Expect F_{n;k;a_1, \ldots, a_k}(t/n)^*F_{n;m;b_1, \ldots, b_m}(t/n) = \\ =
        \sum_{k=0}^n\sum_{1 \leq a_1 < \ldots < a_k \leq n} \Expect F_{n;k;a_1, \ldots, a_k}(t/n)^*F_{n;k;a_1, \ldots, a_k}(t/n),        
    \end{multline*}
    where the last equality follows by Lemmas~\ref{l:adjoint} (Adjoince), \ref{l:super-op} (Super-operator), and~\ref{l:indep} (Independence).
    Indeed, if tuples \((a_1,\ldots,a_k)\) and \((b_1,\ldots,b_m)\) differ, then there exist some <<outer>> composition $B$ (corresponding to the matching tails \((a_{k-s}, \ldots, a_k)\) and \((b_{m-s}, \ldots, b_m)\)), and the <<inner>> composition $C$ which is independent of $B$ since all $\Delta_n(t)$ are independent for pairwise distinct $n \in \N$. Since \(\Expect \Delta_n(t) = 0\) and by Lemma~\ref{l:adjoint} (Adjoince) \(\Expect \Delta_n^*(t) = 0\), we have \(\Expect C = 0\) by Lemma~\ref{l:indep} (Independence) and we may apply Lemma~\ref{l:super-op} (Super-operator) to \(B^*CB\).

    After this the remaining part of the proof just repeats the analogous part from \cite{DzhenzherSakbaev24}.
\end{proof}

\printbibliography

\end{document}